\documentclass[12pt]{article}

\usepackage[a4paper]{geometry}
\usepackage{ajc-FR}
\usepackage{amsmath,amsfonts,latexsym}
\usepackage{graphicx}
\usepackage{placeins}

\Volume{XY(Z)}
\Year{2020}
\firstpage{1}
\lastpage{4}
\Revised{24 May 2020}
\runninghead{\copyright \, F. Rowley 2020\,/\,CWSP1 DRAFT arXiv version, pages 1 -- 4}
\numberwithin{equation}{section}

\newtheorem{theorem}{Theorem}[section]

\newenvironment*{proof}
{\begin{list}{}{\setlength{\leftmargin}{0em}\setlength{\rightmargin}{0em}}
\item[] {\sc Proof:}} {\hfill$\Box$
\end{list}}

\begin{document}

\title{\bf A construction for\\weak Schur partitions}
\author{Fred Rowley\thanks{formerly of Lincoln College, Oxford, UK.}}
\Addr{West Pennant Hills, \\NSW, Australia.\\{\tt fred.rowley@ozemail.com.au}}

% \date{\dateline{submission date}{acceptance date}\\
% \small Mathematics Subject Classifications: comma separated list of
% MSC codes available from http://www.ams.org/mathscinet/freeTools.html}

\date{\dateline{24 May 2020}{DD Mmm CCYY}}

\maketitle

% Papers must include an abstract. The abstract should consist of a
% succinct statement of background followed by a listing of the
% principal new results that are to be found in the paper. The abstract
% should be informative, clear, and as complete as possible. Phrases
% like "we investigate..." or "we study..." should be kept to a minimum
% in favour of "we prove that..."  or "we show that...".  Do not
% include equation numbers, unexpanded citations (such as "[23]"), or
% any other references to things in the paper that are not defined in
% the abstract. The abstract will be distributed without the rest of the
% paper so it must be entirely self-contained.

\begin{abstract}
In 1952, J.H.Braun claimed to have established a formula giving a lower bound for certain partitions of sets of integers into weakly sum-free classes.  However, no proof or supporting construction was published at that time.  

In today's terminology, that claim was equivalent to giving a formulaic lower bound for the weak Schur number $WS(s)$. $WS(s)$ is the maximum number such that there exists a weak Schur partition of the integers from $1$ to $WS(s)$, into $s$ subsets. In a weak Schur partition of a set of integers, there can be no three \textbf{distinct} members $a$, $b$ and $c$ in any subset, such that $a+b=c$. 

An iterative construction described in this paper results in a similar formulaic lower bound.  Although different from that given by Braun, it reproduces the result $WS(6) \ge 554$ implied by his formula, and exceeds it for all larger values of $s$. Various starting points can be used as a basis for the iterations.  

This result itself is no longer remarkable: it has been proven elsewhere, using another construction, that $WS(6) \ge 642$.  Even so, it is hoped that the formula and its underlying construction may nevertheless be of interest to those interested in weak Schur partitions and/or the closely-related linear Ramsey graphs.  

% keywords are optional
\end{abstract}
\bigskip
\small DRAFT \copyright Fred Rowley  May 2020.
\bigskip

\section{Introduction}
This paper originated from attempts to create an iterative construction of weak Schur partitions, which might then be shown to underlie the formula of J.H.~Braun, published in 1952, in \cite{Walker}.  Interestingly, an exact replication was not achieved -- although the result is satisfactory, equalling one of Braun's results, and surpassing his formula as the number of subsets in the partition increases.  

The construction set out below has some elegant aspects, but overall it is not effective in producing very large examples of weak Schur partitions, relative to more modern methods, including those of \cite{EMRS} which do not depend on strict repetition.  However, it may shed some faint light on Braun's claim.

Notation is defined in section 2.  

In section 3, it is proved that, starting from a single specific graph, a series of graphs can be constructed and proven inductively to underpin a formula for $WS(s)$ applicable for all $s \ge 3$.  

In section 4, some very brief conclusions are drawn. 

\section{Definitions and Notation}

In this paper: 

If the set $S$ of integers $\{1, 2, 3, \dots, n\}$ can be partitioned into $s$ disjoint non-empty subsets $S_i$ for $i = 1, 2, 3, \dots, s$, where no subset contains three \textbf{distinct} integers $a, b, c$, such that $a+b=c$, then each such subset is \emph{weakly sum-free} and that partition is a \emph{weak Schur partition}.  Such a partition $p$ may be denoted by $p(s; n)$, or where the order is indeterminate or obvious, simply written as $p(s)$.  The order of the set $S$ is $n$: it is also referred to as the \emph{order} of the partition, and may be written as $\mid p(s) \mid$.  

For any $s$, $WS(s)$ is the maximum value of $n$ such that a weak Schur partition $p(s; n)$ exists.  $WS(s)$ is known as the \emph{weak Schur number}, and its existence is established by Ramsey's Theorem. 

\section{Iterative Construction of Weak Schur Partitions}

\begin{theorem}
  \label{Thm:C-Thm1}
(Iterative Construction Theorem)\\
There is a weak Schur partition of the set of integers $\{1, 2, \dots, 21\}$ into $3$ subsets.  Starting from that partition, it is possible to construct an infinite sequence of weak Schur partitions $p(s)$, with $\mid p(s+1) \mid \, = \, \mid 3p(s) \mid - 1 $.  
\end{theorem}

The theorem depends on a fairly simple construction process and an inductive proof.   

\begin{proof}

The conditions for the induction are as follows:  

\emph{Condition 1:} The partition $p(s)$ is a weak Schur partition with subsets $S_i$ for $i = 1, 2, 3, \dots, s$.

\emph{Condition 2:} There is no pair $a, 2a$, both members of any subset $S_i$ such that $a > 4$. 

\emph{Condition 3:} $p(s)$ has the two additional special properties, that $S_1 \cup \{ \mid p(s) \mid + 2 \}$ is also weakly sum-free, and that $\mid p(s) \mid$ is not a member of $S_1$.  

Assuming these conditions are true for $p(s)$, we can now define the subsets $S'_i$ of the partition $p(s+1)$ by reference to the subsets $S_i$. 

To ease notation, write $m = \mid p(s) \mid$ and $m' = \mid p(s+1) \mid$.

Firstly, let $S'_1 = S_1 \cup \{ m + 2 \} \cup \{ 2m + 2 \} \cup \{ 3m + 4 - a \mid 4 < a \in S_1 \}$.

Then for each other $S_i$, for $2 \le i \le s$, let $S'_i = S_i \cup \{ 3m + 4 - a \mid 4 < a \in S_i \}$.

And lastly let $S'_{s+1} = \{ m + j \mid j = 1, 3, 4, \dots, (m + 1), (m + 3) \}$.

It is simple to verify that this partition is well-defined, and is of order $3m-1$.

We now prove there can be no 'forbidden sum' in any of the subsets $S'_i$ of $p(s+1)$ - that is, no sums $a+b=c$, where either $a \neq b$ or $a = b > 4$.  This will establish that Condition $1$ and Condition $2$ are both true for $p(s+1)$.  

We start with $S'_1$.  

Firstly, we note that the absolute difference between integers $(m+2), (2m+2) \in S'_1$ is equal to $m \notin S'_1$, so that any forbidden sum in the new partition can involve at most one of them.  However, there is no $a \in S'_1$ such that $a+(m+2) \in S'_1$, and therefore there is no $b \in S'_1$ such that $b-(m+2) \in S'_1$. Moreover, there is no $a \in S'_1$ such that $a+(2m+2) \in S'_1$.  If there were, then by the (slightly incomplete) symmetry of the construction, $(3m+4)-(a+2m+2) = (m+2)-a \in S'_1$, which we know is false, by Condition 3.  It follows that there can be no $b \in S'_1$ such that $b-(2m+2) \in S'_1$. We have thus proved that neither of $(m+2), (2m+2)$ is involved in a forbidden sum.

Suppose now that there is a forbidden sum in any $S'_i$, for $i \le s$, then there are two cases.  

In the first case, assume it is of the form $a+b=c$ where $c > 2m+3$.  Then we know that $a \neq b$, and one of those ($b$, say) is also greater than $2m+3$.  If we take the complement of both $b$ and $c$ with respect to $3m+4$, then $4 < 3m+4-b \in S_i$ and $4 < 3m+4-c \in S_i$.  However, the absolute difference between these complements is still $a$, and since both are greater than $4$, neither of them can be equal to twice the other.  This is a contradiction of Condition 1.

In the second case, if $c \le 2m+3$, then we know that $c \le m$.  If $a+b=c$, by Conditions 1 and 2, any sum must be of the form $a+a = c$ with $a \le 4$. 

It easy to verify that $S'_{s+1}$ is (strongly) sum-free.  Thus we have proved that Conditions 1 and 2 are true for $p(s+1)$.

Now note that $3m-1 \in S'_2$, i.e. $3m-1 \notin S'_1$.  Since $m'+2 = 3m+1$, it only remains to prove that $S'_1 \cup \{ 3m+1 \}$ is also weakly sum-free.

Again we proceed to a contradiction.  If there existed $a, b \in S'_1$ such that $a+b=3m+1$, then as before, one of them ($b$, say) would be greater than $2m+1$.  If so, then taking the complement of $b$ and $3m+1$ with respect to $3m+4$ would produce two numbers with absolute difference $a$, both members of $S_1$, and both greater than $4$.  This final contradiction proves the induction, for all $s \ge 3$.  

Now we consider a specific partition $p(3)$, for which we write:

 $p(3; 21) = \{1, 2, 4, 8, 18 \} \cup \{ 3, 5, 6, 7, 19, 20, 21 \} \cup \{ 9, 10, 11, 12, 13, 14, 15, 16, 17 \}$.

It is easy to demonstrate that this is a weak Schur partition, noting that in the third set listed above there are no $a, b, c$, such that $a+b=c$, even where $a$ and $b$ are \textbf{not} distinct.  

We note also that $\{ 1, 2, 4, 8, 18 \} \cup \{ \mid p(3) \mid + 2 \} = \{ 1, 2, 4, 8, 18, 23 \}$ is also a weakly sum-free set, and that $\mid p(3) \mid = 21$ is not a member of $S_1$.  Therefore this $p(3; 21)$ satisfies all the inductive conditions, and this observation completes the proof.

\end{proof}

The recurrence relationship implied by this theorem is that $\mid p(s+1) \mid \, = 3\mid p(s) \mid -\, 1$.

That formula implies $WS(4) \ge 62$, $WS(5) \ge 185$, $WS(6) \ge 554$, and $WS(7) \ge 1661$.  Examination of the Braun formula indicates that it will produce inferior lower bounds on $WS(s)$ for all values of $s > 6$.  

The weak partition into $6$ sets mentioned in \cite{EMRS} is of order $572$, and exceeds the current result.  Both of them exceed the strong partition of order $536$ achieved by Fredricksen and Sweet, mentioned in \cite{FrSw}. Disappointingly, perhaps, the strong partition into $7$ sets, also mentioned in \cite{FrSw}, has order $1680$, exceeding the order of the weak partition produced by this construction. 

However, it is also noted that in a later paper\cite{FR-NLBWSP}, this author has produced a partition of order $642$ in $6$ colours, and $2146$ in $7$ colours, using methods similar to those in \cite{FR-GLRGC}.  By the same method, one can establish the existence of a series of weak Schur partitions with known order, for all higher values of $s$.  

\section{Conclusions}
	
The results of this analysis are of largely historical interest, but nevertheless indicate the extent of the achievement of J.H. Braun in producing his remarkable formula in 1952.  It seems quite likely that he would have used an iterative approach such as that described above, although one cannot be sure, and the precise basis for his result remains a mystery.

% \FloatBarrier
% Use \FloatBarrier if using package "placeins"' - this commmand fixes position of Tables - they cannot float below here.

%\medskip
%%%%%%%%%%%%%%%%%%%%%%%%%%%%%%%%%%%%%%%%%%%%%%%%%%%%%%%
% \subsection*{Acknowledgements}
% I record my warmest thanks .

%%%%%%%%%%%%%%%%%%%%%%%%%%%%%%%%%%%%%%%%%%%%%%%%%%%%%%%
% \bibliographystyle{plain} 
% \bibliography{myBibFile} 

\begin{thebibliography}{5}

\bibitem{EMRS} S.~Eliahou, J.M.~Mar\'{i}n, M.P.~Revuelta, M.I.~Sanz, \newblock Weak Schur numbers and the search for G.W.~Walker's lost partitions, \emph{Comput. Math. Appl.} \textbf{63} (2012), 175-182.

\bibitem{FrSw} H.~Fredricksen and M.M.~Sweet, \newblock Symmetric Sum-Free Partitions and
Lower Bounds for Schur Numbers, \newblock {\em Electron. J. Combin.}, {\bf 7} (2000),  {\#}R32.

\bibitem{FR-GLRGC} F.~Rowley, \newblock A generalised linear Ramsey graph construction, \\
 arXiv preprint arXiv:1912.01164, 2019 - arxiv.org

\bibitem{FR-NLBWSP} F.~Rowley, \newblock New lower bounds for weak Schur partitions, \\
 arXiv preprint arXiv:2011.11292, 2020 - arxiv.org

\bibitem{Walker} G.W.~Walker, A Problem in Partitioning, \emph{Amer. Math. Monthly} \textbf{59} (1952), 253.

\end{thebibliography}
% If you use BibTeX to create a bibliography
% then copy and past the contents of your .bbl file into your .tex file

\end{document}